\input amssym.def 
\input amssym
\magnification=1200
\parindent0pt
\hsize=16 true cm
\baselineskip=13  pt plus .2pt
$ $

\centerline {\bf  On equivariant embeddings of hyperbolic surfaces}
\centerline {\bf into hyperbolic 3-manifolds}

\bigskip

\centerline {Bruno P. Zimmermann}

\bigskip \bigskip

{\bf Abstract.}  We consider the problem of when a closed hyperbolic
surface admits a totally geodesic embedding into a closed hyperbolic
3-manifold, and in particular equivariant versions of such embeddings.
In a previous paper we considered orientation-preserving actions 
on orientable surfaces; in the present paper, we consider large
orientation-reversing actions on orientable surfaces, and also large 
actions on nonorientable surfaces.

\medskip

e-mail:   zimmer@units.it

\bigskip \bigskip

{\bf 1. Introduction}

\medskip

Let $\cal S$ be a closed, orientable or nonorientable  
hyperbolic surface and  $G$ a finite group of isometries of $\cal S$; we say that the
pair $({\cal S},G)$ {\it embeds geodesically} 
if there exists a totally geodesic embedding of 
$\cal S$ into a closed hyperbolic 3-manifold ${\cal M}$ such that the action of $G$ on
${\cal S}$ extends to an isometric action of $G$ on ${\cal M}$. 
Moreover, the pair
$({\cal S},G)$ {\it bounds geometrically} if ${\cal S}$ is the unique totally geodesic 
boundary component of a compact hyperbolic 3-manifold ${\cal M}$ and the action of
$G$ on ${\cal S}$ extends to an isometric action of $G$ on ${\cal M}$.

\medskip

It is shown in [Z3], [GZ] (see also the survey [Z1])
that various infinite series of Hurwitz-actions of maximal possible order $84(g-1)$ on
orientable surfaces of genus $g$ bound geometrically, but that not all Hurwitz
actions bound geometrically. For example, the smallest Hurwitz action of 
${\rm PSL}(2,\Bbb Z_7)$, of order 168 on Klein's quartic of genus 3 does not bound any
compact 3-manifold with a single boundary component. However by [Z2, Corollary 1], every
Hurwitz action embeds geodesically.  Other interesting surfaces which bound
geometrically, with large characterizing finite group actions, are
considered in [F].

\medskip

In [Z2] we considered finite orientation-preserving group actions on orientable
surfaces, in the present paper we consider finite group actions of maximal
possible order on nonorientable surfaces and maximal orientation-reversing actions on
orientable surfaces.  By the formula of Riemann-Hurwitz, the maximal order of a group
of isometries of  a closed orientable hyperbolic surface of genus $g \ge 2$ is
$168(g-1)$; in the following, we call such a group an {\it
orientation-reversing Hurwitz group}  (the orientation-preserving subgroup
of index 2 is then a classical Hurwitz group).  Our first result is the following.

\bigskip
\vfill \eject

{\bf Theorem 1.}   {\sl  Every orientation-reversing Hurwitz
action $({\cal S},G)$ of maximal possible order $168(g-1)$ on a closed orientable
hyperbolic surface $\cal S$ of genus $g$ embeds geodesically into a
closed orientable hyperbolic 3-manifold $({\cal M}, G)$ with an isometric action of
$G$.}

\bigskip

The smallest orientation-reversing Hurwitz action is again on Klein's quartic, for an
action of ${\rm PGL}(2,\Bbb Z_7)$ of order 336.  As noted before,
the smallest Hurwitz action of ${\rm PSL}(2,\Bbb Z_7)$ on Klein's quartic does not
bound. Forgetting about the group action, it remains open whether Klein's quartic bounds
geometrically a hyperbolic 3-manifold; however, it bounds geometrically a hyperbolic
{\it 3-orbifold}.

\bigskip

{\bf Corollary 2.}   {\sl  The surface of every orientation-reversing Hurwitz action
bounds geometrically a compact orientable hyperbolic 3-orbifold. In particular, Klein's
quartic bounds geometrically a compact orientable hyperbolic 3-orbifold.}

\bigskip

We consider nonorientable surfaces now. Again by the formula of
Riemann-Hurwitz, the maximal order of a finite group $\bar G$ of isometries of a
nonorientable hyperbolic surface $\bar {\cal S}$ of genus (or crosscap number) $\bar g
\ge 3$ is $84(\bar g - 2)$, and  we call  such an action $({\bar {\cal S}},\bar G)$  a
{\it nonorientable Hurwitz action} in the following.

\medskip

Let $\cal S$ be the orientable double cover of
$\bar {\cal S}$, a surface of genus $g =
\bar g -1$. The group $\bar G$ lifts to a group 
$G$ of isometries of $\cal S$ of order $168(g-1)$. The orientation-reversing Hurwitz
group $G$ has a normal  subgroup $\Bbb Z_2 = \; <\tau>$ generated by the covering
involution $\tau$ and is isomorphic to $\bar G \times \Bbb Z_2$
where $\bar G$ denotes the orientation-preserving subgroup of index 2 of $G$ now 
(isomorphic to the isometry group $\bar G$ of $\bar {\cal S}$). In particular, $\bar G$
is also a classical Hurwitz group of order $84(g-1)$.

\medskip

By Theorem 1 and its proof,   $({\cal S},G)$
embeds geodesically into a closed orientable  hyperbolic 3-manifold 
$({\cal M}, G)$ with an isometric $G$-action, and $\bar {\cal M} = {\cal
M}/\tau$ is a closed nonorientable hyperbolic 3-manifold with an isometric action of
$\bar G$. Also, 
$({\cal S},G)$ projects to a geodesic embedding of  $({\bar {\cal S}},\bar G)$ in 
$(\bar {\cal M}, \bar G)$, so Theorem 1 implies
(for the second part, see the proofs of 
Corollary 4  or Theorem 2):

\bigskip

{\bf Corollary 3.}   {\sl  Every nonorientable  
Hurwitz action  $({\bar {\cal S}},\bar G)$ on a nonorietable surface ${\bar {\cal S}}$ 
has a totally geodesic, two-sided embedding into a closed nonorientable hyperbolic
3-manifold $(\bar {\cal M},
\bar G)$ with an isometric action of $G$. It has also a totally geodesic, 
one-sided embedding into a closed orientable hyperbolic 3-manifold with an isometric
$G$-action.}

\bigskip

Theorem 1 implies also the following:

\bigskip

{\bf Corollary 4.}   {\sl  For a nonorientable Hurwitz action 
$({\bar {\cal S}},\bar G)$, the associated Hurwitz action
$({\cal S}, \bar G)$ on the orientable double cover of
$\bar {\cal S}$ bounds geometrically an orientable  hyperbolic  
3-manifold $({\cal M'}, \bar G)$ with an isometric action of $\bar G$. }

\bigskip

For a list of the finite groups of order less than $2 \times 10^6$ realizing the maximal
orders for orientation-reversing actions on orientable surfaces, and also for action on
nonorientable surfaces, see [C1, p. 27 and 28]. 

\medskip

In section 3 we generalize our results for larger classes of orientation-reversing
and nonorientable group actions.

\bigskip

{\bf 2. Proof of Theorem 1}

\medskip

Let $G$ be a group of isometries of an orientable hyperbolic surface $\cal S$ of genus
$g$ of maximal possible order $168(g-1)$.  The lift of  $G$ to the universal
covering $\Bbb H^2$ of $\cal S$ is an extended triangle group $[2,3,7]$ generated
by the reflections $r_1, r_2$ and
$r_3$ in the sides of a hyperbolic triangle with angles $\pi/2, \pi/3$ and $\pi/7$,
with a presentation

$$[2,3,7] =  \;\;  < r_1, r_2, r_3 \;|\;  r_1^2 = r_2^2 = r_3^2 = 1, \;  (r_1r_2)^2 =
(r_2r_3)^3 = (r_3r_1)^7 = 1 \;>$$

and a surjection $\phi: [2,3,7]  \to G$ with torsionfree kernel isomorphic to
$\pi_1(\cal S)$ (the universal covering group of $\cal S$).

\medskip

We consider a tetrahedron $\cal T$ as a cone over an
equilateral triangle in the plane and represent it by its  projection to the plane, 
subdividing the equilateral triangle into three triangles which have the
central point in common (the projection of the cone point). In cyclic anti-clockwise
order, we associate labels $r_1, r_2$ and $r_3$ to these three triangles, labels 2, 3
and 7 to the three edges meeting in the cone point, and labels $x$ (a positive integer
$\ge 2$ to be determined), 7 and 7 to the three edges of the (subdivided) equilateral
triangle (such that the edges of the $r_1$-triangle have labels 2, 7 and 7, etc.).

\medskip

A label $y$ at an edge of the tetrahedron stands for a dihedral angle $\pi/y$; 
by Andreev's theorem on hyperbolic polyhedra ([A], [V], [T], [RHD]), $\cal T$ can
be represented by a hyperbolic tetrahedron ${\cal T}_x$ in $\Bbb H^3$ with orthogonally
truncated vertices of types [2,3,7] (the central vertex), [$x$,3,7],
[7,7,7], and a fourth vertex of type [$x$,2,7] which is truncated only if $x > 2$ 
(spherical if $x = 2$).

\medskip

Let $r_1, r_2$ and $r_3$ denote the reflections in the faces of ${\cal T}_x$ with these
labels, and $r_4$ the reflection in the fourth face of ${\cal T}_x$. The Coxeter group
$C_x$ generated by the reflections in the four faces of ${\cal T}_x$ has a presentation

$$C_x \; = \; < r_1, r_2, r_3, r_4 \; | \; r_1^2 = r_2^2 = r_3^2 = r_4^2 = 1, \; 
(r_1r_2)^2 = (r_2r_3)^3 = (r_3r_1)^7 = 1,$$  
$$ (r_1r_4)^7 = (r_2r_4)^x = (r_3r_4)^7 = 1  \;>$$ 

(cf [B] for a presentation of the tetrahedral group associated to ${\cal T}_x$, the
orientation-preserving subgroup of index 2 of the Coxeter group; for Poincar\'e's
theorem on fundamental polyhedra, its history and references, see [R]).

\medskip

The extended triangle group [2,3,7] is a subgroup of $C_x$ acting on the
orthogonally truncating hyperbolic plane of the cone point of  ${\cal T}_x$. 
We want to extend the surjection 
$\phi: [2,3,7] \to G$ to a surjection $\psi: C_x  \to G$ whose kernel $K$ acts freely on
$\Bbb H^3$. Setting $\psi(r_i) = \phi(r_i)$, $1 \le i \le 3$, we still have
to  define $\psi(r_4)$.

\bigskip

The subgroup of [2,3,7] generated by
$r_1$ and $r_3$ is a dihedral group of order 14 whose image under $\phi$  in
$G$ we denote by $\Bbb D_{14}$. The dihedral group $\Bbb D_{14}$ consists of 7
rotations and 7 reflections, and we choose for $\psi(r_4)$ a reflection in $\Bbb
D_{14}$  different  from  $\phi(r_1), \phi(r_2)$ (which is maybe not even in $\Bbb
D_{14}$),  and  $\phi(r_3)$. Then $\psi (r_1r_4)$ and 
$\psi (r_3r_4)$ have order 7, and we define $x$ as the order   
of $\psi (r_2r_4)$. With these choices of  $\psi(r_4)$ and $x$, the surjection 
$\psi: C_x \to G$ extends the surjection $\phi: [2,3,7] \to G$ such that the kernel $K$
of $\psi$ acts freely on $\Bbb H^3$.

\medskip

We denote by $\bar \Bbb H^3$ the subset of $\Bbb H^3$ obtained by
truncating $\Bbb H^3$ by the truncating planes of the vertices of the tetrahedron
${\cal T}_x$ and  their images under the action of $C_x$.  Then $\bar \Bbb H^3$ is
invariant under the actions of $C_x$ and $K$,  and ${\cal M}' = \bar \Bbb H^3/K$ is a
compact orientable hyperbolic 3-manifold with totally geodesic boundary and an action of
$G \cong C_x/K$, and the boundary component of ${\cal M}'$ corresponding to the cone
point of ${\cal T}_x$ is the surface $\cal S$ with the action of $G$ defined by $\phi$. 
The double of ${\cal M}'$ and its $G$-action along the
boundary is a closed orientable hyperbolic 3-manifold $\cal M$ with a totally geodesic
embedding of $({\cal S},G)$.

\bigskip

{\it Proof of Corollary 2.}  Let  $({\cal S},G)$  be an orientation-reversing Hurwitz
action of maximal possible order $168(g-1)$ of a closed orientable hyperbolic surface
$\cal S$. By Theorem 1 and its proof,  $({\cal S},G)$ embeds geodesically into a closed
orientable hyperbolic 3-manifold $({\cal M}, G)$ with an isometric $G$-action; 
moreover, by splitting $\cal M$ along $\cal S$ one obtains  a compact hyperblic
3-manifold with an isometric $G$-action whose boundary consists of two totally geodesic
copies of $\cal S$ invariant under the action of $G$. Choosing one ot these two copies,
the
$G$-action on it contains an orientation-reversing involution $\tau$. Identifying
points on this copy by $\tau$,  one obtains an orientable hyperbolic 3-orbifold with
exactly one totally geodesic boundary component isometric to $\cal S$ (a hyperbolic
3-manifold if $\tau$ has no fixed points).

\bigskip

{\it Proof of Corollary 4.}  Let $G$ denote the lift of $\bar G$ to the
orientable double cover $\cal S$ of $\bar{\cal S}$, with covering involution 
$\tau \in   G =  \bar G \times \Bbb Z_2$ generating $\Bbb Z_2$.  Since, by
Theorem 1,  $({\cal S},G)$ embeds geodesically into a closed orientable hyperbolic
3-manifold  $({\cal M},G)$, one can  split $\cal M$ along  $\cal S$
and then close up one of the two arising boundary components, identifying points on
it by the orientation-reversing fixed-point free involution
$\tau \in G$. The result is an orientable  hyperbolic 3-manifold $\cal M'$ with
an isometric action of ${\bar G} \subset G$ and a unique, 
totally geodesic boundary component $\cal S$.

\bigskip \bigskip
\vfill  \eject

{\bf 3. A generalization}

\medskip

We denote by $[p,q,r]$ the extended triangle group generated by the reflection in the
sides of a hyperbolic triangle with angles $\pi/p, \pi/q$ and $\pi/r$, and by $(p,q,r)$
its orientation-preserving subgroup of index 2.  Let 
$\phi: [p,q,r] \to  G$  be a surjection with torisonfree kernel $K$ onto a finite group
$G$; then  ${\cal S} = \Bbb H^2/K$ is a closed surface with an isometric action of $G$.

\bigskip

{\bf Theorem 2.}   {\sl  Let $({\cal S},G)$ be an isometric action of a finite group
$G$ on a closed surface ${\cal S} = \Bbb H^2/K$,  associated to a surjection $\phi:
[p,q,r]  \to  G$ with torisonfree kernel $K$. 

\smallskip

i) If $\cal S$ is orientable then $({\cal S},G)$ embeds geodesically into a closed
orientable hyperbolic 3-manifold $({\cal M},G)$ with an isometric $G$-action.

\smallskip

ii) If $\cal S$ is nonorientable then $({\cal S},G)$ has a totally geodesic, two-sided
embedding into a closed nonorientable hyperbolic 3-manifold $({\cal M},G)$ with an
isometric action of $G$. It has also a totally geodesic, one-sided embedding into a
closed orientable hyperbolic 3-manifold with an isometric action of $G$.}

\bigskip

We note that, by the lists in [C2], most of the large orientation-reversing group
actions on orientable surfaces and actions on nonorientable surfaces, up
to genus 300, are of the type of Theorem 2.

\bigskip

{\it Proof of Theorem 2.}   The proof is analogous to the proof of Theorem 1; 
instead of the cone over a triangle (the tetrahedron $\cal T$) we consider the
double cone over a triangle now, a polyhedron $\cal P$ with six faces and five vertices.
As in the proof of Theorem 1, we associate labels $r_1, r_2$ and $r_3$ to the three faces
around the upper cone point and labels $r_1', r_2'$ and $r_3'$ to the three faces around
the lower cone point, in such a way that $r_1$ and $r_2'$ have a common edge with label
$p$,  $r_2$ and $r_3'$ a common edge with label $q$ and 
$r_3$ and $r_1'$ a common edge with label $r$.  The three edges around the upper
cone point have labels $p, q$ and $r$ (starting with the edge betwwen $r_1$ and $r_2$),
and the same for the three edges around the lower cone point (starting with the edge
between $r_1'$ and $r_2'$).  The two cone points of the polyhedron $\cal P$ are then of
type $[p,q,r]$, the three remaining vertices of types $[p,p, q,q]$,  $[q,q,r,r]$ and 
$[r,r,p,p]$.

\medskip

By Andreev's theorem, the polydedron $\cal P$ with these labels $p, q$ and $r$ can be
realized by a hyperbolic polyhedron with dihedral angles $\pi/p,  \pi/q$ and $\pi/r$   
and with orthogonally truncated vertices. Let $C$ be the Coseter group
generated by the reflections (also denoted by $r_i$ and $r_i'$) in the six faces of the
polyhedron $\cal P$ . At the
truncating plane of the upper cone point, there acts the extended triangle group
$[p,q,r]$ generated by $r_1, r_2$ and $r_3$, and we extend the 
surjection $\phi: [p,q,r] \to G$ to a surjection  $\psi: C \to G$ by setting 
$\psi(r_i')  = \psi (r_i) = \phi(r_i)$, $\; 1\le i \le 3$.

\medskip

The proof  of part i) of Theorem 2 
finishes now exactly as the proof of Theorem 1 (and also of the first statement of ii).

\medskip

For the proof of ii), suppose that $({\cal S},G)$ is an
action on a nonorientable surface $\cal S$. Let $\tilde {\cal S} \to {\cal S}$ be the
orientable double cover of $\cal S$, with covering involution $\tau$.  The action of
$G$ on $\cal S$ lifts to an action of $\tilde G = G \times \Bbb Z_2$ on 
$\tilde {\cal S}$ where $\tau$ generates $\Bbb Z_2$ and $G$ is interpreted as
the orientation-preserving subgroup of $\tilde G$ now. By i),  
$(\tilde {\cal S}, \tilde G)$ has a totally geodesic embedding into a closed orientable
3-manifold  $(\tilde {\cal M}, \tilde G)$ such that $\tau \in \tilde G$ acts freely on 
$\tilde {\cal M}$ (also, $({\cal M}, G) = (\tilde {\cal M}/\tau, \; \tilde
G/\tau)$  is a nonorientable closed hyperbolic 3-manifold with a totally geodesic
embedding of $({\cal S}, G) = (\tilde {\cal S}/\tau, \tilde G/\tau)$, reproving the
first statement of ii).

\medskip

By splitting $\tilde {\cal M}$ along $\tilde {\cal S}$, we obtain an orientable
hyperbolic 3-manifold with an isometric action of $\tilde G$ whose totally geodesic
boundary consists of two copie of 
$\tilde {\cal S}$. Identification of  points on these two boundary components by the
fixed point free, orientation-reversing involution $\tau$ gives a closed orientable
hyperbolic 3-manifold with an isometric action of $G$ and  two totally geodesic
embeddings of $({\cal S}, G) = (\tilde {\cal S}/\tau, \tilde G/\tau)$.

\bigskip \bigskip

\centerline {\bf References}

\bigskip

\item {[A]}  E.M. Andreev,  {\it On convex polyhedra in Lobachevsky space.} 
Math. USSR Sb. 10  (1970),  413-440

\item {[B]}  L.A. Best,  {\it  On torsion-free discrete subgroups of PSL(2, \Bbb C)
with compact orbit space.} Can. J. Math.  23  (1971), 451-460

\item {[C1]}  M. Conder,  {\it  Maximal automorphism groups of symmetric Riemann
surfaces with small genus.}  J. Algebra 114  (1988),  16-28

\item {[C2]}  M. Conder,  {\it  Large group actions on surfaces.} 
Contemp. Math. 629  (2014),  77-97

\item {[F]}  L. Ferrari,  {\it Geometric bordisms of the Accola-Maclachlan, Kulkarni
and Wiman type II surfaces.}  arXiv:2312.01821

\item {[GZ]}  M. Gradolato, B. Zimmermann,  {\it  Extending finite group actions on
surfaces to hyperbolic 3-manifolds.}  Math. Proc. Camb. Phil. Soc. 117  (1995),
137-151

\item {[R]}  J.G. Ratcliffe,  {\it  Foundations of Hyperbolic Manifolds.}  Graduate
Texts in Mathematics 149, Springer-Verlag 1994

\item {[RHD]}  R.K.W. Roeder, J.H. Hubbard, W.D. Dunbar,  {\it Andreev's theorem on
hyperbolic polyhedra.}  arXiv:math/0601146

\item {[T]} W. Thurston,  {\it The geometry and topology of 3-manifolds.} 
Lecture Notes, Dept. of Math. Princeton Univ. 1977

\item {[V]} E.B. Vinberg,  {\it Geometry II.}  Encyclopaedia Math. Sciences vol. 29,
Springer 1993

\item {[Z1]} B. Zimmermann,  {\it Large finite group actions on surfaces: Hurwitz
groups, maximal reducible and maximal handlebody groups, bounding and non-bounding
actions.}  

arXiv:2110.11050

\item {[Z2]}  B. Zimmermann, {\it  On geodesic embeddings of hyperbolic surfaces into
hyperbolic 3-mani-folds.}  arXiv:2401.06651

\item {[Z3]}  B. Zimmermann, {\it  Hurwitz groups and finite group actions on
hyperbolic 3-manifolds.}  

J. London Math. Soc. 52  (1995), 199-208

\bye